\documentclass[12pt] {article}
\usepackage{amssymb,amsxtra,amsmath,amsthm}
\def \ZZ{{\mathbb{Z}}}

\def \RR{{\mathbb{R}}}

\def \FF{{\mathbb{F}}}
\def \PP{{\mathbb{P}}}

\def \CCC{{\mathcal{C}}}
\def \OOO{{\mathcal{O}}}
\def \TTT{{\mathcal{T}}}

\def \mm{{\mathfrak m}}

\def \RRRR{{\mathfrak R}}

\def\det{\mathop{\mathrm{det}}}

\def\mod{\mathop{\mathrm{mod}}}

\def\vert{\mathop{\mathrm{vert}}}
\def\edge{\mathop{\mathrm{edge}}}
\def\diag{\mathop{\mathrm{diag}}}
\def\star{\mathop{\mbox{\Large$*$}}}
\def\Star{\mathop{\mbox{\LARGE$*$}}}
\def\card{\mathop{\mathrm{card}}}

\def\diag{\mathop{\mathrm{diag}}}

\def\dim{\mathop{\mathrm{dim}}}

\def\Cl{\mathop{\mathrm{Cl}}}
\def\val{\mathop{\mathrm{val}}}

\def\Gamma{\varGamma}
\def\Delta{\varDelta}

\begin{document}

\begin{center}
{\Large {\bf The stabilizers in a Drinfeld modular group of the vertices of its Bruhat-Tits tree: an elementary approach}}\\
\bigskip
{\tiny {\rm BY}}\\
\bigskip
{\sc A. W. Mason}\\
\bigskip
{\small {\it Department of Mathematics, University of Glasgow\\
Glasgow G12 8QW, Scotland, U.K.\\
e-mail: awm@maths.gla.ac.uk}}\\
\bigskip
{\tiny {\rm AND}}\\
\bigskip
{\sc Andreas Schweizer\footnote{During the writing of this paper
the second author was working at Academia Sinica in Taipei,
supported by grant 99-2115-M-001-011-MY2 from the
National Science Council (NSC) of Taiwan.}}\\
\bigskip
{\small {\it Department of Mathematics,\\
Korea Advanced Institute of Science and Technology (KAIST),\\
Daejeon 305-701, SOUTH KOREA\\
e-mail: schweizer@kaist.ac.kr}}
\end{center}
\begin{abstract}
Let $K$ be an algebraic function field of one variable with constant
field $k$ and let $\CCC$ be the Dedekind domain consisting of all
those elements of $K$ which are integral outside a fixed place
$\infty$ of $K$.  When $k$ is {\it finite} the group $GL_2(\CCC)$
plays a central role in the theory of Drinfeld modular curves analagous
to that played by $SL_2(\ZZ)$ in the classical theory of modular forms.
When $k$ is {\it finite} (resp. {\it infinite}) we refer to a group
$GL_2(\CCC)$ as an {\it arithmetic} (resp. {\it non-arithmetic})
{\it Drinfeld modular group}. Associated with $GL_2(\CCC)$ is its
{\it Bruhat-Tits tree}, $\TTT$. The structure of the group is
derived from that of the quotient graph $GL_2(\CCC)\backslash \TTT$.
Using an elementary approach which refers explicitly to matrices we determine the structure of all the vertex
stabilizers of $\TTT$. This extends results of Serre, Takahashi and
the authors. We also determine all possible valencies of the
vertices of $GL_2(\CCC)\backslash \TTT$ for the important special
case where $\infty$ has degree $1$.
\end{abstract}
\newpage

\begin{center}{\bf \large Introduction}
\end{center}

\noindent Let $K$ be an algebraic function field of one variable
with constant field $k$. As usual we assume that $k$ is algebraically
closed in $K$. Let $\infty$ be a fixed place of $K$ of degree
$\delta$ and let $\CCC$ be the set of all those elements of $K$
which are integral outside $\infty$. Then $\CCC$ is a Dedekind
domain whose unit group, $\CCC^*$, is $k^*$. (For the simplest
example let $K=k(t)$, the rational function field over $k$. Then,
when $\infty$ corresponds to the usual ``point at infinity" of
$k(t)$, $\CCC$ is the polynomial ring $k[t]$.) Here our focus of
attention is the group $G=GL_2(\CCC)$. When $k$ is {\it finite} (the
{\it arithmetic} case) this group plays a central role [1], [2] in
the theory of {\it Drinfeld modular curves}, analogous to that of
the classical modular group $SL_2(\ZZ)$ in the theory of modular
forms. When $k$ is finite we refer to $GL_2(\CCC)$ as an ({\it
arithmetic}) {\it Drinfeld modular group}. Otherwise we call
$GL_2(\CCC)$ a {\it non-arithmetic Drinfeld modular group}.
\\ \\
\noindent Let $K_{\infty}$ be the completion of $K$ with respect to
$\infty$. The group $GL_2(K_{\infty})$ acts [9, Chapter II, Section
1.1] on its associated {\it Bruhat-Tits building} which in this case
is a {\it tree}, $\TTT$. Classical Bass-Serre theory [9, Theorem
10, p.39] shows how to obtain a presentation for $GL_2(\CCC)$ from
the structure of the quotient graph $GL_2(\CCC)\backslash \TTT$.
Serre's approach [9, 2.1, p.96] to the structure of $\TTT$ is based on the theory of vector bundles.
 Here, as in [3], we adopt a more elementary approach which explicitly refers to matrices. One of our aims in so doing
 is to
make our results more accessible to group theorists less familiar with algebraic geometry.
Moreover the approach based on vector bundles does not provide more detailed versions of our principal results
(which describe the group-theoretic structure of the vertex stabilizers in $GL_2(\CCC)$).  Serre [9, Theorem 9, p.106] has determined the ``shape" of
$GL_2(\CCC)\backslash \TTT$. He has shown that this quotient is the
union of a central subgraph $X$, of bounded width, together with a
number of (pairwise disjoint) infinite half-lines (or {\it rays}).
(``Bounded width" refers to the geodesic length in
$GL_2(\CCC)\backslash \TTT$.) In the arithmetic case $X$ is finite
and there are only finitely many rays. The structure of the rays (by
which we mean the structure of the stabilizers of vertices and edges
of $\TTT$ which project onto those of the rays) is well understood
[9, p.118]. However much less is known about $X$, except for a
number of special cases. For example the structure of $X$ is known
[4], when $g=0$ and $K \cong k(t)$, and [9, 2.4.3, p.114], [6],
when $g=0$, $K \not\cong k(t)$ and $\delta=2$. In addition Takahashi
[12] has determined $X$ precisely when $g=1$ and $ \delta=1$. Among
the few known [9, p.97] general properties of $X$ is that it
contains at least one vertex whose stabilizer is isomorphic to
$GL_2(k)$.
\\ \\
\noindent A number of important properties of $GL_2(\CCC)$ derive
from $X$. For example, $GL_2(\CCC)$ has a free quotient whose rank
is that of $\pi_1(X)$, the fundamental group of $X$. (See [9, p.43].)
The theory of Drinfeld modular curves [2] provides a formula [5] for
the (finite) rank of $\pi_1(X)$, for the case where $k$ is finite.
\\ \\
\noindent The principal aim in this paper is to determine the structure
of the stabilizers of all the vertices of $X$ thus extending the
above results of Serre, Takahashi and the authors. Our elementary approach shows that this structure depends entirely on the
nature of the eigenvalues of the relevant matrices. Our main result
is the following. Suppose that there exists a quadratic polynomial
over $k$ which has no roots in $k$.
(Then $k$ is said to be not {\it quadratically closed}.)
\\ \\
\noindent {\bf Theorem A.} \it Let $G_v$ be the stabilizer in $G$ of
$v \in \vert(\TTT)$.
\\ \\
\noindent (a) Suppose that the
eigenvalues of every matrix in $G_v$ lie in $k$ and that at least
one such matrix has distinct eigenvalues. Then
$$ G_v /N \cong k^* \times k^*,$$
where  $N \cong V^+$, the additive group of $V$, a
finite-dimensional vector space over $k$.
\\ \\
\noindent (b) Suppose that every matrix in $G_v$ has repeated
eigenvalues in $k$. Then $$G_v = Z \times N,$$ where $N$ is as
above and $Z \cong k^*$.
\\ \\
\noindent (c) Suppose that $G_v$ contains a matrix whose eigenvalues
are distinct and do not lie in $k$. Then $G_v$ is isomorphic to one
of the following.
\begin{itemize}
\item[(i)] $L^*$, where $L$ is a separable quadratic extension of $k$,
\item[(ii)] $GL_2(k)$,
\item[(iii)] $Q^*$, the units of a quaternion division algebra $Q$ over
$k$.
\end{itemize}

\noindent (d) Suppose that every matrix of $G_v$ has repeated
eigenvalues and that one such matrix has eigenvalues not in $k$.
Then $\mathrm{char}(k)=2$ and
$$G_v \cong F^*,$$
\noindent where $F$ is a totally inseparable finite extension of $k$ of
degree $2^m$, where $m \geq 1$. Moreover $2^{m-1}$ divides the degree of
every place of $K$.
\\ \\
\noindent \rm Non-arithmetic groups are more complicated than their
arithmetic counterparts. For example, stabilizers of types (c)(iii)
and (d) {\it only} arise when $k$ is infinite. The situation when
$k$ is {\it quadratically closed} turns out to be much simpler.
\\ \\
\noindent \rm Our results show that although edge stabilizers are of the same type they are subject to some restrictions. For the remainder of the paper we focus on the valency
(or degree) of the projection of a vertex of $\TTT$ in
$GL_2(\CCC)\backslash\TTT$. For one important special
case we determine all the possibilities.
\\ \\
 \noindent {\bf Theorem B.} \it Suppose that $\delta=1$. Let
 $\widetilde{v}$ denote the projection of a vertex $v$ of $\TTT$ in
 $GL_2(\CCC)\backslash\TTT$. Then
 $$\val(\widetilde{v})=1,\;2,\;3\;\ or\;\card(k)+1.$$

\noindent \rm Of particular interest here are the {\it isolated} vertices
in the quotient graph for the following reason. If a vertex $v$ and edge
$e$ project onto such a vertex and its only incident edge in
$GL_2(\CCC)\backslash\TTT$, then from standard Bass-Serre theory
[9, Theorem 13, p.55] it follows that
$$ GL_2(\CCC) \cong H \star_{\quad L }K,$$
\noindent where $H=G_v$ and $L=G_e$, the stabilizers of $v$ and $e$
in $\TTT$, with $H \neq L$ and $K \neq L$. For isolated vertices the
situation is again less complicated for arithmetic groups. For
example we prove here that, in contrast to the non-arithmetic case,
isolated vertices exist in the quotient graph {\it only} when
$\delta=1$.
\\ \\
In the case where $k$ is finite, $GL_2(\CCC)\backslash\TTT$ is a
particularly interesting object from a number theoretic point of
view, because it encodes a lot of information concerning the
Drinfeld modular curve associated with $GL_2(\CCC)$.
\par
In [7] we work out a precise connection between the elliptic
points of this Drinfeld modular curve (that is, points on the
Drinfeld upper half-plane with nontrivial stabilizer under the action
of $G$) and vertices of $G\backslash\TTT$ with certain stabilizers.
In the case $\delta=1$ we also get a precise relation between the
elliptic points and the isolated vertices of $G\backslash\TTT$. This
in turn then yields information on the possibilities for
$PGL_2(\CCC)$ to decompose as a free product.
\par
In contrast to the mainly group theoretic approach in the present
paper, the arguments used in [7] are largely number theoretic,
involving properties of ideal class groups in function fields and
$L$-functions.
\\ \\

\begin{center}{\bf \large 1. Preliminaries}
\end{center}
 \par\noindent For convenience we list at this point the notation
which will be used throughout:

\begin{tabular}{ll}
$k$             & a field;\\
$\FF_q$         & the finite field of order $q$; \\
$K$             & an algebraic function field of one variable with constant field $k$;\\
$g=g(K)$        & the genus of $K$;\\
$\infty$        & a chosen place of $K$;\\
$\delta$        & the degree of the place $\infty$;\\
$\nu$           & the discrete valuation of $K$ defined by $\infty$;\\
$\pi$           & a local parameter at $\infty$ in $K$;\\
$\OOO$          & the valuation ring of $\infty$ in $K$;\\
$\mm=(\pi)$     & the maximal ideal of $\OOO$;\\
$k_{\infty}$     &  the residue field, $\OOO/\mm$;\\
$K_{\infty}$     & the completion of $K$ with respect to $\infty$;\\
$\OOO_{\infty}$  & the completion of $\OOO$ with respect to $\infty$;\\
$\TTT$          & the Bruhat-Tits tree of $GL_2(K_\infty)$;\\
$\CCC$          & the ring of all elements of $K$ that are integral outside $\infty$;\\
$G$             & the group $GL_2(\CCC)$;\\
$G_w$           & the stabilizer in $G$ of $w \in \vert(\TTT)\cup\edge(\TTT)$;\\
$\widetilde{w}$ & the image in $G\backslash \TTT$ of $w \in \vert(\TTT)\cup\edge(\TTT)$;\\
$Z$             &  the set of scalar matrices in $G$;\\
$Z_{\infty}$     &  the set of scalar matrices in $GL_2(K_{\infty})$.\\
\end{tabular}
\\ \\
\noindent Our basic reference for the theory of algebraic function
fields is Stichtenoth's book [11]. We recall that $k_{\infty}$ is a
finite extension of $k$ of degree $\delta$. Associated with the
group $GL_2(K_{\infty})$ is its {\it Bruhat-Tits building} which in
this case is a {\it tree}, $\TTT$. See [9, Chapter II, Section 1].
It is known [9, Corollary, p.75] that $G$ acts on $\TTT$ {\it
without inversion}. Classical Bass-Serre theory [9, Theorem 13,
p.55] shows how the structure of $G$ can be derived from that of the
quotient graph $G\backslash
\TTT$. The ``shape" of this graph is described in the following.
\\ \\
\noindent {\bf Theorem 1.1.(Serre)} {\it To each element $\sigma$
of $\Cl(\CCC)$, the ideal class group of $\CCC$, there corresponds a
ray (i.e. an infinite half-line without backtracking), $R(\sigma)$, with
terminal vertex $v_\sigma$, together with a subgraph $X$ such that
$$ G \backslash \TTT=\left(\bigcup_{\sigma \in \Cl(\CCC)}R(\sigma)\right)\cup
X,$$ where \begin{itemize}
\item[(i)] $X$ is bounded (with respect to geodesic length in
$G\backslash \TTT$),
\item[(ii)] $R(\sigma) \cap R(\tau)=\emptyset \;(\sigma \neq \tau)$,
\item[(iii)] $\vert(X) \cap \vert(R(\sigma))=\{v_{\sigma}\}$,
\item[(iv)] $\edge(X)\cap\edge(R(\sigma))=\emptyset.$
\end{itemize}}

\noindent Serre's approach [9, Theorem 9, p.106] uses the theory of
vector bundles. For a more elementary approach which refers
specifically to matrices see [3, Theorem 4.7]. A presentation for
$G$ can be derived from a {\it lift}
$$j: \TTT_0 \longrightarrow \TTT,$$
where $\TTT_0$ is a maximal tree of $G\backslash \TTT$. Each ray
$R(\sigma)$ is then realised as a subgraph
of $\TTT$.\\
\noindent  The structure of the stabilizers of the vertices of
$\TTT$ which project onto those of each $R(\sigma)$ is well
understood. See, for example, [3, Theorem 4.2]. However much less
is known about $X$. One of the few known [9, p.97] properties of
$X$ is that it contains at least one vertex whose stabilizer is
isomorphic to $GL_2(k)$. The precise structure of $X$ is only known
for a number of special cases, for example, when $g=0$ [4], [6]
and when $g=1$ and $ \delta=1$ [12].
\\ \\
\noindent A number of
important properties of $G$ are derived from $X$. For example, it is
known [9, p.43] that $G$ has a free quotient whose rank is that of,
$\pi_1(X)$, the fundamental group of $X$. The theory of Drinfeld
modular curves [2] provides a formula [5] for the (finite) rank of
$\pi_1(X)$, for the case where $k$ is finite. When $k$ is infinite this
rank is only known in a number of particular instances [4],[6],[12].
\\ \\
\noindent We require a detailed model for $\TTT$. For convenience we
make use of that used by Takahashi [12]. The {\it vertices} of $\TTT$
are the left cosets of $Z_{\infty}GL_2(\OOO_{\infty})$ in
$GL_2(K_{\infty})$, where $Z_{\infty}GL_2(\OOO_{\infty})$ is the
subgroup generated by $GL_2(\OOO_{\infty})$ and $Z_{\infty}$. Recall
that
$$\OOO=\{ z \in K:\nu(z) \geq 0\},$$
$\nu(\pi)=1$ and that $\OOO/\mm\cong k_{\infty}$. In addition $\nu(c)
\leq 0$, for all $c \in \CCC$ and $\OOO\cap\CCC=k$. The {\it edges}
of $\TTT$ are defined in the following way. The vertices
$g_1Z_{\infty}GL_2(\OOO_{\infty})$ and
$g_2Z_{\infty}GL_2(\OOO_{\infty})$ are {\it adjacent} if and only if
$$g_2^{-1}g_1 \equiv  \left[\begin{array}{lll} \pi & z\\[10pt]
0 & 1
\end{array}\right]\; \mathrm{or}\: \left[\begin{array}{lll} \pi^{-1} & 0\\[10pt]
0 & 1
\end{array}\right]\; (\mod \;Z_{\infty}GL_2(\OOO_{\infty})).
$$
for some $z \in k_{\infty}$. Then $\TTT$ is a tree on which $G$ acts
(by left multiplication) {\it without inversion}. It is clear that
$\TTT$ is {\it regular} in the sense that the edges of $\TTT$
incident with each $v \in \vert(\TTT)$ are in one-one correspondence
with the elements of $\PP_1(k_{\infty})$. It is also clear that if $e
\in\edge(\TTT)$ has vertex extremities $u,v$
then $G_e=G_u\cap G_v$.\\
\noindent The model used by Serre [9, Chapter II, 1.1] is different
but equivalent. A {\it lattice} $L$ is a $\OOO_{\infty}$-submodule
of $K_{\infty}^2$ of rank $2$ and the {\it class} containing $L$ is
the set $\{zL:\; z\in K_{\infty}^*\}$. Then $G$ acts (naturally) on
the set of lattice classes which for this model are the vertices of
$\TTT$. Let $gZ_{\infty}GL_2(\OOO_{\infty})$ be any left coset and
let $$ g=\left[\begin{array}{lll} a & b\\[10pt]
c & d
\end{array}\right].$$
\noindent If $\Lambda_g$ denotes the lattice class containing the
lattice generated by $(a,c)$ and $(b,d)$ then the correspondence
$$gZ_{\infty}GL_2(\OOO_{\infty})\longleftrightarrow \Lambda_g$$
demonstates that the models for $\TTT$ are equivalent.\\ \\
\noindent Since $\OOO_{\infty}$ is a PID we may represent every
coset of $Z_{\infty}GL_2(\OOO_{\infty})$ by an element of the form
$$\left[\begin{array}{lll} \pi^n & z\\[10pt]
0 & 1
\end{array}\right],
$$
for some $ n \in \ZZ$ and $z \in K_{\infty}$. We denote this vertex
of $\TTT$ by $v(n,z)$. It is clear that
$$ v(n,z)=v(m,z') \Longleftrightarrow n=m \;\mathrm{and}\; \nu(z-z') \geq n.$$
We may assume therefore that $ z \in K$. \noindent When $v=v(n,z)$
we put $$G_v=G(n,z).$$
\\
\noindent {\bf Definition.} Let $H$ be
any subgroup of $G$. Vertices $u,v$ of $\TTT$ are said to be
$H$-{\it equivalent}, written $$ u \equiv v\;(\mod \;H),$$ \noindent
if and only if $u=h(v)$, for some $h \in H$.
\\ \\
\noindent If, as above, $u,v$ are $H$-{\it equivalent} it is clear
that $G_v=hG_uh^{-1}$ and that $\bar{u}=\widetilde{v}$. A similar
definition can be made for the elements of $\edge(\TTT)$.
\\ \\\noindent
We record the following well-known result.
\\ \\
\noindent {\bf Lemma 1.2.} \it With the above notation, let
$v=v(n,z)$ and
$$ M=\left[\begin{array}{lll} a & b\\[10pt]
c & d
\end{array}\right]\in G.$$
\noindent Then $M \in G(n,z)$ if and only if
\begin{itemize}
\item[(i)]$\nu(c) \geq -n,$
\item[(ii)]$ \nu(a-zc) \geq 0,\;\nu(d+zc)\geq 0,$
\item[(iii)] $\nu(b+z(a-d)-z^2c)\geq n.$
\end{itemize}

\noindent \rm Lemma 1.2 shows that there exists a constant
$\kappa=\kappa(v)$ such that $\nu(x) \geq \kappa$, for every entry
of the matrices in $G_v$. Hence the entries of the matrices in $G_v$
lie in a {\it bounded} subset of $K_{\infty}$ with respect to the
metric on $K_{\infty}$ defined by $\nu$. (Serre [9, Proposition 2,
p.76] proves an alternative version of this result.) For each $n \in
\ZZ$ let $$\CCC(n)=\left\{c \in \CCC:\nu(c) \geq -n\right\}.$$
\noindent By the Riemann-Roch Theorem [11, I.5.15, p.28] $\CCC(n)$
is a {\it finite-dimensional} vector space over $k$. It follows
that, when $k$ is finite, each $G_v$ is {\it finite}. The following notation is useful for our purposes.
\\ \\
\noindent {\bf Notation.} We put
$$M(n,z,\alpha,\beta,c)=\left[\begin{array}{cc} \alpha + cz & b\\[10pt]
c &\beta-cz \end{array}\right],$$ \noindent where \begin{itemize}
\item[(i)] $\nu(c) \geq -n,$
\item[(ii)] $c,cz \in \CCC,$
\item[(iii)] $b=(\beta-\alpha)z-cz^2 \in \CCC,$
\item[(iv)] $\alpha,\beta \in k^*.$
\end{itemize}
\noindent By Lemma 1.2 it is clear that $M(n,z,\alpha,\beta,c) \in
G(n,z)$ and that $\det(M(z,n,\alpha,\beta,c))=\alpha\beta$. Let
$$V=\left\{ c \in \CCC: M(n,z,1,1,c) \in G(n,z)\right\}.$$
\noindent From the above, if $ c \in V$, then $c \in \CCC(n),
\;cz,cz^2 \in \CCC$ and $V$ is a finite-dimensional $k$-vector
space. We note also that $M(n,z,\alpha,\alpha,0)=
\diag(\alpha,\alpha) \in G(n,z)$ and that
$$M(n,z,\alpha,\alpha,c) \in G(n,z)\;\Longleftrightarrow\; c\in \CCC(n), \; cz,cz^2 \in \CCC.$$

\noindent {\bf Lemma 1.3.} \it Let $ v \in \vert(\TTT)$ and let $M
\in G_v$. Then the characteristic polynomial of $M$ has coefficients
 in $k$.
\\ \\
\noindent {\bf Proof.} \rm Let $$ M=\left[\begin{array}{lll} a & b\\[10pt]
c & d
\end{array}\right].$$
\noindent It suffices to prove that ${\rm tr}(M)=a+d \in k$, since
$\det(M) \in k^*$. By Lemma 1.2 $\nu(a+d) \geq 0$ and so $a+d \in
\CCC\cap\OOO=k$.
\hfill $\Box$
\\ \\
\noindent \rm Our study of stabilizers depends crucially on the next
result.
\\ \\
\noindent {\bf Lemma 1.4.} \it Let $w \in
\vert(\TTT)\cup\edge(\TTT)$. Suppose that the matrices $M_1,M_2 \in
G_w$. Then, if
$$\det(\alpha_1M_1+\alpha_2M_2) \in k^*,$$
where $\alpha_1,\alpha_2 \in k$, then
$$\alpha_1M_1+\alpha_2M_2 \in G_w.$$

\noindent {\bf Proof.} \rm We may assume that $w \in \vert(\TTT)$.
Then $M_i$ ``fixes" some coset $gZ_{\infty}GL_2(\OOO_{\infty})$,
say. Hence $g^{-1}M_ig\in Z_{\infty}GL_2(\OOO_{\infty})$ and so
$g^{-1}M_ig \in GL_2(\OOO_{\infty})$, since $\det(M_i) \in k^*$
($i=1,2$). It follows that
$$g^{-1}(\alpha_1M_1+\alpha_2M_2)g \in M_2(\OOO_{\infty}).$$
\hfill $\Box$
\\
\noindent We note that $Z=\{\alpha I_2: \alpha \in k^*\}$ and
consequently that
$$Z\leq G_w,$$
\noindent for all
$w \in \vert(\TTT)\cup\edge(\TTT)$.
\\ \\
\noindent Let the characteristic polynomial $\chi_M(X)$ of $M \in
GL_2(K_{\infty})$ have coefficients in $k$. It is clear that, for
all $\alpha,\beta \in k$, with $\beta \neq 0$,
$$\det(\alpha I_2+\beta M)=\beta^{2} \chi_M(-\alpha \beta^{-1}) \in k.$$
{\bf Notation.} We put
$$I(M)=\{ M' =\alpha I_2+\beta M: \alpha,\beta \in k,\;\det(M')\in k^*\}.$$
\noindent It is clear by Lemmas 1.3 and 1.4 that, if $M \in G_v$ then
$$Z \leq I(M) \leq G_v.$$ This will be the starting point for our study of
$G_v$. We record an immediate consequence of Lemma 1.4.
\\ \\
\noindent {\bf Lemma 1.5.} \it Let $M \in G$ and
$w\in\vert(\TTT)\cup\edge(\TTT)$. Then
$$I(M) \cap G_w=Z\;\mathrm{or}\;I(M).$$

\noindent {\bf Lemma 1.6.} \it Suppose that $M \in G_v \backslash
Z$. Then there are three possibilities.
\begin{itemize}

\item[(i)] If $M$ has distinct eigenvalues in $k$, then
$$I(M) \cong k^* \times k^*.$$
\item[(ii)] If $M$ has repeated eigenvalues in $k$, then
$$I(M) \cong k^* \times k^+,$$
\noindent where $k^+$ is the additive group of $k$.
\item[(iii)] If $M$ has an eigenvalue $\lambda \notin k$, then
$$ I(M) \cong k(\lambda)^*.$$
\end{itemize}

\noindent {\bf Proof.} \rm For part (i) $M$ is conjugate to
$D_0=\diag(\lambda_1,\lambda_2)$, where $\lambda_1,\lambda_2 \in
k^*$ , with $\lambda_1 \neq \lambda_2$. Then $$ I(M) \cong I(D_0) =
\{ \diag(\lambda,\mu): \lambda,\mu \in k^*\}.$$

\noindent Part (ii) is very similar. Here $$I(M) \cong Z \times B,$$
where $$B =\left\{\left[\begin{array}{lll} 1 & x\\[10pt]
0 & 1
\end{array}\right]: x \in k\right\}.$$
\noindent For part (iii) we note that $\det (\alpha I_2+\beta M)=0$
only when $\alpha=\beta=0$. It is clear then that in this case
$I(M)$ is the multiplicative group of a quadratic extension of $k$.\hfill $\Box$\\ \\
\noindent We record separately a special case of Lemma 1.6.
\\ \\
\noindent {\bf Lemma 1.7.} \it Suppose that $k=\FF_q$. Then with the
above notation
$$|\;I(M)|=\left\{\begin{array} {cll}(q^2-1)&,&the\; eigenvalues\; of\; M\;are\;not\; in\; \FF_q \\
(q-1)^2&,&the\; eigenvalues\;
of\;M\;are\;distict\;and\;in\;\FF_q\\q(q-1)&,&otherwise\end{array}
\right. $$
\noindent \rm Note that $I(M) \cong \FF_{q^2}^*$, when $|I(M)|=q^2-1.$
\\ \\
\noindent We conclude this section by stating a well-known result
which follows immediately from the definition of $\TTT$.
Let $\RRRR\subset K$ denote a complete set of
representatives for $\OOO/\mm\;(\cong k_{\infty})$.
\\ \\
\noindent
{\bf Lemma 1.8.} \it The vertices in $\TTT$ adjacent to $v(n,z)$ are
$$v(n-1,z)\;\; and \;\;v(n+1,z+u\pi^n),$$
where $u \in \RRRR$.
\\ \\

\begin{center}{\bf \large 2. Vertex stabilizers}
\end{center}

\noindent \rm This section is devoted to our principal aim, namely
the determination of the structure of each $G_v$. We require the 
following.
\\ \\
\noindent {\bf Definition.} The field $k$ is called
{\it quadratically closed} if and only if every quadratic
polynomial over $k$ has a root in $k$.
\\ \\
\noindent Every algebraically closed field is, of course,
quadratically closed. Examples of quadratically closed fields which
are not algebraically closed include the separable closure of an
imperfect field of odd characteristic, and
$$\bigcup_{n\geq\;0}\FF_{q(n)},$$ where $q(n)=q^{2^n}$. Examples of
fields which are not quadratically closed include all subfields of
$\RR$, all finite fields, and all fields with a discrete valuation,
in particular local fields, algebraic number fields, and algebraic
function fields over any constant field.
\\ \\
\noindent {\bf Theorem 2.1.} \it  Suppose that all the matrices in
$G_v$ have eigenvalues in $k$ with at least one with distinct
eigenvalues. Then there are two possibilities.
\begin{itemize}
\item[(a)] There exists a normal subgroup $N$ of $G_v$, such that
\begin{itemize}
\item[(i)]$N \cong V^+$, the additive group of a finite-dimensional
$k$-vector space $V$ ,
\item[(ii)] $G_v/N\cong k^*\times k^*$.
\end{itemize}
\noindent In this case there is a homomorphism $\theta : k^*\times
k^* \rightarrow \mathrm{Aut}(V^+)$ given by
$$\theta((\alpha,\beta)): v \mapsto (\alpha\beta^{-1})v \;\;(v \in
V).$$
\item[(b)] Only when $k$ is quadratically closed, $$G_v \cong GL_2(k).$$
\end{itemize}

 \noindent {\bf Proof.}  \rm  Let $M_0 \in G_v$ have distinct eigenvalues
in $\lambda, \mu \in k^*$.
Then we may replace $G_v$ with a conjugate $G'_v$ (over
$GL_2(K)$) containing
$$I(D_0) = D,$$
where $D_0=\diag(\lambda,\mu)$ and
$D=\{\diag(\sigma,\tau):\sigma,\tau \in k^*\}$. \noindent Now
$G'_v$ satisfies the same hypotheses as $G_v$. In addition
Lemmas 1.3, 1.4 and 1.6 apply to $G'_v$. Let
$$ M=\left[\begin{array}{lll} a & b\\[10pt]
c & d\end{array}\right]\in G'_v.$$ \noindent By multiplying $M$
and $\diag(\lambda,1)$ and considering traces we conclude by Lemma
1.3 that
$$\lambda a+d\in k$$
for all $\lambda \in k^*$. It follows that $a,d \in k$ and hence
that $bc \in k$. (Note that the hypotheses ensure that $k \neq
\FF_2$.) There are two possibilities.\\ \\ \noindent {\bf (a)}
Suppose that $bc=0$, for all $M \in G'_v$. It follows then that
either $c=0$, for all $M \in G'_v$, or $b=0$, for all $M \in
G'_v$. We will assume the former.

\noindent We define a subset $V$ of $K$ by $$ x \in V
\Longleftrightarrow \left[\begin{array}{lll}1 & x\\[10pt]
0 & 1\end{array}\right] \in G'_v.$$ \noindent Let
$$B=\left\{\left[\begin{array}{lll}1 & x\\[10pt]
0 & 1\end{array}\right]:x\in V\right\}.$$ \noindent It is clear that
 $$ G'_v \left\{\left[\begin{array}{lll}\alpha & x\\[10pt]
0 & \beta\end{array}\right]:\alpha,\beta\in k^*,\;x\in V\right\}.$$
\noindent By Lemma 1.4 $V$ is a vector space over $k$. Moreover $B
\cong V^+$. \noindent Under the conjugacy between $G_v$ and
$G'_v$ $B$ is conjugate to a subgroup $B'$, say, of $G_v$. We
now apply Lemma 1.2 to the entries of $B'$. It follows that there
exists $ t \in K^*$ and a real constant $\kappa$ such that, for all
$x \in V$,
\begin{itemize}
\item[{(i)}]  $tx \in \CCC$,
\item[{(ii)}] $\nu(tx) \geq \kappa.$
\end{itemize}
\noindent Then $V$ is a
{\it finite-dimensional} $k$-space by the Riemann-Roch theorem. \\ \\
\noindent {\bf (b)} Suppose that there exists $$ M_0=\left[\begin{array}{lll} a_0 & b_0\\[10pt]
c_0 & d_o\end{array}\right]\in G'_v$$ for which $b_0,c_0 \in
k^*$. Then, for any $M \in G'_v$, as above, since $MM_0 \in
G'_v$ it follows from the above that $bc_0, cb_0 \in k$. We
deduce that
$$ G'_v=\left\{\left[\begin{array}{cc}\alpha & \beta b_0\\[10pt]
\gamma c_0& \delta\end{array}\right]:\alpha,\beta, \gamma, \delta
\in k,\; \alpha\delta-\beta\gamma \neq 0 \right\}.$$ \noindent By
multiplying $M_0$ with a suitable element of $D$ we may assume that
$b_0c_0=1$. It is clear then that $$G'_v \cong GL_2(k),$$
\noindent in which case $k$ is quadratically closed.
\hfill $\Box$
\\ \\
\noindent The arithmetic version of Theorem 2.1 is as follows.
\\ \\
\noindent {\bf Corollary 2.2.} \it Suppose that $k= \FF_q$. Then,
with the hypotheses of Theorem 2.1,
$$G_v / N \cong \FF_q^*\times \FF_q^*,$$
\noindent where
$$N \cong V^+,$$ the additive group of a finite-dimensional
$\FF_q$-vector space $V$ (and hence elementary $p$-abelian).
\noindent In addition
$$|\;G_v|=(q-1)^2q^n ,$$ where $n=\dim_q(V) \geq 0.$
\\ \\
\noindent \rm Theorem 2.1 applies in particular [3, Theorem 4.2] to
the vertices of $\TTT$ which map onto those of the rays $R(\sigma)$
(in Serre's Theorem). If $v_1,v_2$ map onto adjacent vertices of
$R(\sigma)$ in $G\backslash \TTT$, then $G_{v_i} \cong V_i^+
\rtimes(k^* \times k^*)\; (i=1,2)$. Moreover (relabelling if
necessary) it is known that $G_{v_1} \leq G_{v_2}$ and that
$\dim(V_2)-\dim(V_1)=\delta$. Consequently $\dim(V)$ in Theorem 2.1
can be arbitrarily large.
\\ \\
\noindent We record a known ``minimal
size" example of a stabilizer as described by Theorem 2.1. Suppose
that $g=\delta=1$. Then $\CCC$ is the coordinate ring of an {\it
elliptic} curve. In which case there exist $x,y \in \CCC$, where
$\nu(x)=-2$ and $\nu(y)=-3$ which satisfy a {\it Weierstrass
equation}, $F(x,y)=0$, for which
$$\CCC=k[x,y].$$\noindent Takahashi [12] has completely determined
$G \backslash \TTT$ in this case. In particular [12, Theorem 5] he
has shown that
$$G(2,\pi^{-1}+ \pi\lambda) \cong k^* \times k^*,$$ \noindent
whenever there exist $\mu_1,\mu_2 \in k$, with $\mu_1 \neq \mu_2$,
such that $$F(\lambda,\mu_1)=F(\lambda,\mu_2)=0.$$

\noindent {\bf Theorem 2.3.} \it Suppose that every matrix in $G_v$
has repeated eigenvalues in $k$. Then
$$ G_v = Z\times B,$$
where $ B \cong V^+$, the additive group of a finite-dimensional
vector space $V$ over $k$.
\\ \\
\noindent {\bf Proof.} \rm We recall that $Z \leq G_v$. We may
assume that $G_v \neq Z$. Then $G_v$ contains a non-central matrix.
As in the proof of Theorem 2.1 we replace $G_v$ with a conjugate
$G'_v$ which contains $Z$ and a matrix
$$ \left[\begin{array}{lll} \alpha & x\\[10pt]
0 & \alpha\end{array}\right],$$ where $\alpha \in k^*$ and $x \neq
0$. Now as above $G'_v$ satisfies the same hypotheses as $G_v$
and Lemmas 1.3, 1.4 and 1.6 apply to it. By Lemma 1.4 it follows
that, for all $\lambda \in k^*$,
$$ U=\left[\begin{array}{lll} 1 & \lambda x\\[10pt]
0 & 1\end{array}\right] \in G'_v.$$

\noindent Let $$ M=\left[\begin{array}{lll} a & b\\[10pt]
c & d\end{array}\right]\in G'_v.$$ \noindent Then there are two
possibilities.
\\ \\
\noindent Suppose that $\mathrm{char}(k) \neq 2.$ Since the eigenvalues
are repeated, we have $(a+d)^2 =4(ad-bc)$, that is $(a-d)^2=-4bc$. But
$UM \in G'_v$ and so
$$(a-d+\lambda cx)^2=-4c(b+\lambda dx).$$
\noindent It follows that
$$2c(a-d)+c^2\lambda x=-4cd,$$
and hence that $c=0$.
\\ \\
\noindent Suppose that $\mathrm{char}(k)=2.$ The hypotheses ensure
that in this case $a=d$, for all $M$. Again by considering $UM$ it
follows that $c=0$.
\\ \\
\noindent In all cases therefore $c=0$ and
hence $a=d$. Defining $V$ and $B$ as in the
proof of Theorem 2.1 we deduce that
$$G'_v= \left\{\left[\begin{array}{lll} \alpha & b\\[10pt]
0 & \alpha\end{array}\right]: \alpha \in k^*,\;b \in V \right\}.$$
\noindent Again $V$ is a $k$- vector space which is
{\it finite-dimensional} by the Riemann-Roch theorem.
\hfill $\Box$
\\ \\
\noindent The arithmetic version of Theorem 2.3 is as follows.
\\ \\
\noindent {\bf Corollary 2.4.} \it Suppose that $k=\FF_q$. Then,
with the hypotheses of Theorem 2.3,
\begin{itemize}
\item[(i)]$G_v\cong \FF_q^* \times V^+,$ where $V^+$ is the additive
group of a vector space over $\FF_q$, \item[(ii)]$|\;G_v|=(q-1)q^n$,
where $n=\dim(V)\geq 0$. \end{itemize}

\noindent \rm The images in $G\backslash \TTT$ of the vertices to
which Theorem 2.3 and Corollary 2.4 apply lie in $X$ (in Serre's
theorem), except when
$k=\FF_2$.
\\ \\
\noindent We record `minimal size" examples of stabilizers as
described in Theorem 2.3. Suppose that $g >0$ and that $\delta=1$.
Then by the {\it Weierstrass gap theorem} [11, I.6.7, p.32]  $1$ is
a {\it gap number} for $\nu$, i.e. there is {\it no} $x \in \CCC$
for which $\nu(x)=-1$. It follows readily from Lemma 1.2 that
\begin{itemize}
\item[(i)]  $$G(1,\pi^{-1})\cong k^*,$$
\item[(ii)] $$G{(0,\pi^{-1})}\cong k^* \times k^+.$$
\end{itemize}

\noindent \rm The remaining cases are more complicated. Let $\overline{K_\infty}$
denote the algebraic closure of $K_{\infty}$. The group
$GL_2(\overline{K_\infty})$ acts as a group of linear fractional transformations
on $\PP^1(\overline{K_\infty})=\overline{K_\infty} \cup \{\infty\}$ in the usual way. For
each subgroup $H$ of $GL_2(\overline{K_\infty})$ and each $z \in \PP^1(\overline{K_\infty})$
let $H_z$ denote the stabilizer in $H$ of $z$.
We record without proof the following.
\\ \\
\noindent {\bf Lemma 2.5.} \it Let $g \in GL_2(\overline{K_\infty})$. Then
$\left[\begin{array}{c}\alpha\\\beta\end{array}\right]$ is an
eigenvector of $g$ if and only if $$g(\gamma)=\gamma,$$ (as a linear
fractional transformation) where $\gamma=\alpha/\beta$.
($\gamma=\infty$, when $\beta=0$).
\\ \\
\noindent {\bf Theorem 2.6.} \it Suppose that $G_v$ contains a
matrix $M_0$ with eigenvalues $\alpha_0,\beta_0$, where
$\alpha_0,\beta_0 \notin k$ and $\alpha_0 \neq \beta_0$. Let
$L=k(\alpha_0)(=k(\beta_0))$ and let $\sigma$ denote the (Galois)
$k$-automorphism of (the quadratic extension) $L$.
\noindent Then there exist $x,y$, for which either $x=y=0$ or
$xy=\kappa \in k^*$, such that
$$G_v \cong\left\{\left[\begin{array}{lll} \lambda & x\mu\\[10pt]
y\mu^\sigma &\lambda^\sigma \end{array}\right]: \lambda,\mu \in L,\;
\lambda\lambda^\sigma \neq \kappa\mu\mu^\sigma \right\}.$$

\noindent There are then two possibilities.
\begin{itemize}\item[(i)] $$G_v\cong L^*.$$
\item[(ii)]  There exists a quaternion algebra $Q$ over
$k$ such that $$G_v \cong Q^*.$$
\end{itemize}

\noindent {\bf Proof.} \rm The nontrivial $k$-automorphism of $L$
extends naturally to a $K$-automorphism of the (quadratic) extension
$LK/K$. Note that $\beta_0=\alpha_0^\sigma$. Let
$\left[\begin{array}{c}\alpha\\\beta\end{array}\right]$ and
$\left[\begin{array}{c}\alpha^\sigma\\\beta^\sigma\end{array}\right]$
be corresponding eigenvectors of $M_0$. Then there exists $g \in
GL_2(\overline{K_\infty})$ such that
$$(I(M_0))^g=\bar{Z}=
\{\diag(\lambda,\lambda^{\sigma}):\lambda \in L^*\}\leq (G_v)^g.$$
\noindent We may assume that $g$ maps
$\left[\begin{array}{c}\alpha\\\beta\end{array}\right]$ and
$\left[\begin{array}{c}\alpha^\sigma\\\beta^\sigma\end{array}\right]$
onto $\left[\begin{array}{c}1\\0\end{array}\right]$ and
$\left[\begin{array}{c}0\\1\end{array}\right]$, respectively. Let
$X=(G_v)^g$. Although $(G_v)^g$ is not necessarily contained in $G$,
Lemmas 1.3, 1.4 and 1.6 apply to the matrices in $X$. Then by
Lemma 2.5 and the above
$$X_0=X_\infty=\bar{Z}.$$
\noindent We first of all dispose of the simplest possibility, namely
$\bar{Z}=X$, i.e. $G_v=I(M_0) \cong L^*$.
(See Lemma 1.6(iii).) We may suppose then that from now on there exists
$$ M=\left[\begin{array}{lll} a & b\\[10pt]
c & d\end{array}\right]\in X\backslash\bar{Z}.$$
\noindent By considering the trace of the product $DM$, where
$D=\diag(\lambda,\lambda^{\sigma})$, with $\lambda \in L^*$.
It follows that
$$\lambda a+\lambda^{\sigma} d \in k,\;\mathrm{for\; all}\; \lambda \in L.$$
\noindent Now $\lambda(a+d) \in L$. From the case where
$\lambda \neq \lambda^{\sigma}$ we deduce that $a,d \in L$.
\\ \\
\noindent Also $\lambda a+ \lambda^{\sigma} a^{\sigma} \in k$ and so
$$\lambda^{\sigma}( d-a^{\sigma}) \in k,\; \mathrm{for\; all}\; \lambda \in
L.$$ \noindent We conclude that $d=a^{\sigma}$ and hence that $bc\in k$.
\\ \\
\noindent We now show that $bc=0$ if and only if
$b=c=0$. Suppose then that $c=0$. By an obvious modified version of
Lemma 1.4 it follows that
$$ B=\left[\begin{array}{lll} 1 & b\\[10pt]
0 &1\end{array}\right]\in X.$$ \noindent Hence $B \in X_{\infty}$
and so $B \in X_0$, from the above. We deduce that $b=0$. Similarly
$c=0$ when $b=0$.
\\ \\
\noindent We may now suppose that $bc \in k^*$. Again by a modified
version of Lemma 1.4 and the above
$$ N=\left[\begin{array}{lll} 0 & b\\[10pt]
c & 0\end{array}\right]\in X.$$ \noindent Let
$$ N_1=\left[\begin{array}{lll} 0 & b_1\\[10pt]
c_1 & 0\end{array}\right]\in X.$$
\noindent  Now $NN_1 \in X$. It
follows from the above that $bc_1,b_1c \in L$ with
$(bc_1)^{\sigma}=b_1c$ and $bc,b_1c_1 \in k^*$. We deduce that
$$b_1=b\beta \;\mathrm{and}\;c_1=c\beta^{\sigma} ,$$
for some $\beta \in L$, i.e. $N_1=\diag(\beta,\beta^{\sigma})N$.
\\ \\
\noindent Let $D_0=\diag(\alpha_0,\beta_0)$. By Lemma 1.4 and the
above it is clear that $X$ consists of all those matrices of nonzero
determinant in the ($4$-dimensional) $k$-vector space $Q$ spanned by
$I_2,\; D_0,\; N$ and $D_0N$. (Let $b=x$ and $c=y$.) Equivalently
$X$ is the set of units in the $k$-algebra $Q$ generated by
$I_2,D_0$ and $N$. It is clear that $Q$ is central. To prove that
$Q$ is quaternion it suffices therefore to prove that, if $J$ is a
nonzero two-sided ideal in $Q$, then $J=Q$.
\\ \\
\noindent Let
$$ M=\lambda I_2+\mu D_0+\nu N+\rho N_0 \in J,$$
where $N_0=D_0N$ and $\lambda,\mu,\nu,\rho \in k$, with
$(\lambda,\mu,\nu,\rho)\neq (0,0,0,0)$. Then
$$C=MD_0-D_0M=\nu (ND_0-D_0N)+\rho(N_0D_0-D_0N_0)
=(\beta-\alpha)\left[\begin{array}{lll} 0 & u\\[10pt]
v &0\end{array}\right] \in J,$$
where $u=x(\nu +\rho\alpha)$ and $v=-y(\nu+\rho\beta)$. Then
$C^2$ is a scalar matrix. If $C^2$ is nonzero it follows that
$J=Q$. If on the other hand it is zero then $\nu=\rho=0$. Thus
$$\lambda I_2+\mu D_0 \in J,$$
where $(\lambda,\mu)\neq (0,0)$. In this case $M$ is invertible
and again $J=Q$.
\hfill $\Box$
\\ \\
\noindent The proof of Theorem 2.6 is
based on a representation of $G_v$ which derives from a quadratic
extension of $k$ generated by the eigenvalues of a particular
element of $G_v$. In general {\it any} quadratic extension of $k$
arises in this way. For example the stabilizer of $G(0,0)$ is
$GL_2(k)$ (by Lemma 1.2) and  every quadratic extension of $k$
is generated by the eigenvalues of some matrix in $GL_2(k)$.\\
By Theorem 2.6 there are three possibilities for $G_v$.
\\ \\
\noindent {\bf Corollary 2.7.} \it With the notation of Theorem 2.6,
$$G_v \cong L^* \Longleftrightarrow x=y=0.$$

\noindent \rm Such vertices exist, for example, when $g=\delta=1$.
See [12, Theorem 5].
\\ \\
\noindent {\bf Notation.} Let $N_{L/k}:L \longrightarrow k$
denote the {\it norm} map, that is $N_{L/k}(a)=a\sigma(a)$.
\\ \\
\noindent {\bf Corollary 2.8.} \it With the notation of Theorem 2.6,
suppose that $\kappa \in N_{L/k}(L^*)$. Then $$ G_v \cong
GL_2(k).$$ \noindent {\bf Proof.} \rm  With the notation of Theorem
2.6, replacing $N$ with $YN$, for some $Y \in \bar{Z}$, we may
assume that $xy=\kappa=1$. It is easily verified that
$$G_v \cong\left\{\left[\begin{array}{lll} \lambda & x\mu\\[10pt]
y\mu^\sigma &\lambda^\sigma \end{array}\right]: \lambda,\mu \in L,\;
\lambda\lambda^\sigma \neq \mu\mu^\sigma
\right\}\cong\left\{\left[\begin{array}{lll} \lambda & \mu\\[10pt]
\mu^\sigma &\lambda^\sigma \end{array}\right]: \lambda,\mu \in L,\;
\lambda\lambda^\sigma \neq \mu\mu^\sigma \right\}.$$

\noindent The latter group is the set of matrices of nonzero
determinant (in $k$) in the $k$-vector space, $V$, say, spanned by
$\{I_2,D_0,N_1,D_0N_1\}$, where
$$N_1=\left[\begin{array}{lll} 0 & 1\\[10pt]
1 &0 \end{array}\right].$$
\noindent Now let
$$M=\left[\begin{array}{lll} \alpha_0 & \alpha_0^{\sigma}\\[10pt]
1 &1 \end{array}\right].$$ \noindent It is readily verified that
$MD_0M^{-1},\;MN_1M^{-1} \in M_2(k)$. Since the $k$-dimension of $V$
is $4$ it follows that $$MVM^{-1} = M_2(k).$$ \noindent The
result follows.
\hfill $\Box$
\\ \\
\noindent Corollary 2.8 arises in all cases since (to repeat from
the above) the stabilizer of $G(0,0)$ is always $GL_2(k)$. The
remaining case can be stated without proof.
\\ \\
\noindent {\bf Corollary 2.9.} \it With the notation of Theorem 2.6,
suppose that $\kappa \notin N_{L/k}(L)$. Then $G_v$ is
isomorphic to the set of units in the quaternion algebra, $Q$, over
$k$, with $k$-basis $\{I_2,D_0,N,D_0N\}$.
\\ \\\noindent \rm The existence of vertex stabilizers as described in
Corollary 2.9 imposes restrictions on $K$.
\\ \\
\noindent {\bf Lemma 2.10.} \it A necessary condition for the
existence of a non-split quaternion algebra $Q$ for which $Q^*$ is a
vertex stabilizer is that all places of $K$ have even degree. \rm
\\ \\
{\bf Proof.} We start from $Q^*\subseteq GL_2(\CCC)\subseteq
GL_2(K)$, fix a $k$-basis of $Q$ and consider its norm form. After
tensoring up to $K$ the norm form will have a nontrivial zero;
either because the $k$-basis becomes linearly dependent over $K$, or
if not because the above inclusion implies that $Q\otimes_k K\cong
M_2(K)$. Thus for every place $v$ of $K$ the norm form of $Q$ has a
nontrivial zero in the completion $K_v$ and hence also in its
valuation ring $\OOO_v$. Obviously we can achieve that the
coordinates of such a zero are not all in the maximal ideal. So
after reduction the norm form of $Q$ will have a nontrivial zero in
the residue field $k_v$. In other words, $k_v$ is a splitting field
of $Q$. Hence the degree of $k_v$ over $k$ must be even by [8, Lemma
in Section 13.4].
\hfill $\Box$
\\ \\
\noindent Lemma 2.10 explains the absence of stabilizers of this
type from Takahashi's list [12, Theorem 5] (since there
$\delta=1$). We are also able to prove the following interesting
existence theorem.
\\ \\
\noindent  {\bf Theorem 2.11.} \it The set of units of every
non-split quaternion algebra $Q$ over $k$ occurs as a vertex
stabilizer of $GL_2(\CCC)$ for a nonrational genus zero function
field $K$ and a place $\infty$ of degree $2$. \rm
\\ \\
{\bf Proof.}
If the characteristic of $k$ is different from $2$, there exists a $k$-basis
$\{1,i,j,ij\}$ of $Q$ with $i^2 =-\rho$, $j^2 =-\sigma$, and $ij=-ji$.
Since $Q$ is non-split, $-\rho$ is not a square in $k$, and for all
$x,y\in k$ the element $y+xi+j$ has norm $y^2 +\rho x^2 +\sigma\neq0$.
We define
$$\CCC:=k[X,Y]\ \ \hbox{\rm  with}\ \ Y^2 +\rho X^2 +\sigma=0.$$
Then $K=k(X,Y)$ is a non-rational function field of genus 0 (compare
[6, Theorem 1.1 case (i)]). By [6, Lemma 2.5] and [6, Remark
2.7] there exists a vertex of $\TTT$ whose stabilizer in $GL_2(\CCC)$
is isomorphic to $Q$.
\par
If the characteristic is $2$, there exists a $k$-basis $\{1,i,j,ij\}$ of
$Q$ with $i^2 +i=\rho$, $j^2 =\sigma$, and $ij=j(i+1)$. For all $x,y\in k$
the element $y+xi+j$ has norm $y^2+xy+\rho x^2 +\sigma\neq 0$. Moreover,
$\rho$ is not of the form $\alpha^2 +\alpha$ for any $\alpha\in k$. Then
the claim follows from [6, Theorem 1.1 (iv)], [6, Lemma 2.6] and
[6, Remark 2.7].
\hfill $\Box$
\\ \\
\noindent \rm  In particular [6, Example 3.3] for this case, when
$k=\RR$, there is a vertex whose stabilizer is {\it Hamilton's
quaternions}. We record separately the arithmetic version of Theorem
2.6.
\\ \\
\noindent {\bf Corollary 2.12.} \it Suppose that $k=\FF_q$. With the
hypotheses of Theorem 2.6,
$$G_v \cong \FF_{q^2}^*\;\; or \;\; GL_2(\FF_q).$$
\noindent {\bf Proof.} \rm In this case $N_{L/k}$ is always
{\it surjective}. (Alternatively we could argue that there are no
non-split quaternion algebras over a finite field, or that a
function field with finite constant field always has places of odd
degree.) \hfill $\Box$
\\ \\
\noindent There remains one further possibility.
\\ \\
\noindent {\bf Theorem 2.13.} \it Suppose that (a) every matrix in
$G_v$ has repeated eigenvalues and (b) the eigenvalues of at least
one matrix in $G_v$ do not lie in $k$. Then $k$ is infinite,
$\mathrm{char}(k)=2$ and
$$ G_v \cong F^*,$$ where \begin{itemize} \item[(i)] $F$ is a finite,
totally inseparable extension of $k$, \item[(ii)] $u^2 \in k$, for
all $u\in F$,
\item[(iii)] $|F:k|=2^m$, for some $m>0$,
\item[(iv)] $2^{m-1}$ divides the degree of every place of $K$, in
particular $2^{m-1}$ divides $\delta$ and $2g-2$.
\end{itemize}

\noindent {\bf Proof.} \rm The fact that $k$ is infinite of characteristic
$2$ follows from Lemma 1.3. Let
$$M=\left[\begin{array}{lll} a &b\\[10pt] c &d \end{array}\right]\in G_v.$$
Then the hypotheses imply that
\begin{itemize}
\item[(i)] $a=d$,
\item[(ii)] either $b=c=0$ or $bc \neq 0$,
\item[(iii)] $b=\gamma c$, where $\gamma$ is a nonzero constant
determined by $G_v$. \end{itemize} \noindent As usual we denote the
set of all squares of elements of $K$ by $K^2$. Suppose that $\gamma
\in K^2$. Then $$\det(M)=a^2+\gamma c^2 \in K^2\cap k.$$ \noindent
\noindent Since $k$ is algebraically closed in $K$ it follows that
the eigenvalues of $M$ lie in $k$. Hence $$\gamma \notin K^2.$$
\noindent Every matrix in $G_v$ has trace $0$ and hence is a multiple
of its inverse. From this one easily sees that $G_v$ is abelian.
Consequently, if $M\in G_v$ is a matrix with eigenvalues outside $k$,
then $G_v$ must be contained in the commutator of $K(M)$ in $GL_2(K)$,
that is, in $K(M)$.
\par
Let $$M_i=\left[\begin{array}{lll} a_i &b_i\\[10pt]
c_i &d_i \end{array}\right]\;(i=1,2)$$ \noindent Then
$\det(M_1+M_2)=(a_1+a_2)^2+\gamma(c_1+c_2)^2 \in k$. It is clear
that $\det(M_1+M_2)= 0$ if and only if $M_1=M_2$. Now $M_1+M_2 \in
G_v$, when $\det(M_1+M_2) \neq 0$, by Lemma 1.4. It is clear that
$G_v$ is the multiplicative group of an extension $F$ of $k$ which
satisfies property (ii). We use $[a,c]$ as a shorthand for
$$\left[\begin{array}{lll} a &\gamma c\\[10pt]
c &a \end{array}\right]. $$
\noindent Then from the above
\begin{itemize}
\item[(i)] If $[a_i,c_i] \in G_v,\;i=1,2$, then $[a_1+a_2,c_1+c_2]\in G_v$.
\item[(ii)] If $[a,c],[a_0,c] \in G_v$, then $a+a_0 \in k$.
\end{itemize}
\noindent Let
$$V=\left\{ c \in \CCC:\;[*,c] \in G_v\right\}.$$
\noindent It is clear from the above, Lemma 1.2 and
the Riemann-Roch theorem that $V$ is a finite-dimensional vector
space over $k$. It follows that $F$ is a {\it finite} extension of
$k$. Property (ii) implies that $F/k$ is totally inseparable,
whence (i). Consequently $F$ satisfies property (iii).
\\ \\
\noindent Choose $M_0 \in F \backslash k$. Then, for all $M \in F$,
there exist $\lambda,\mu \in K$ such that
$$M= \lambda I_2+\mu M_0.$$
\noindent For part (iv) let $L=KF$. Then $L/K$ is a constant field
extension (with constant field $F$) which by the above is also
quadratic. Some of the unusual properties of such extensions (where
the constant field extension is inseparable) can be found in [10].
Alternatively see
[13, Section 8.6] for the same material in English.
\\ \\
Combining [10, Korollar 9] (or [13, Corollary 8.6.15]) with [10,
Lemma 4] (or [13, Theorem 8.6.8]) it follows that $2^{m-1}(=
|F:k|/|L:K|)$ divides the degree of every place of $K$.
\hfill $\Box$
\\ \\
\noindent We note that, in contrast to Theorem 2.13,  the (quadratic)
extension in Corollary 2.7 (and Corollary 2.12) is {\it separable}.
\\ \\
\noindent It is known that stabilizers of the type
described in Theorem 2.13 do exist. For example, when $g=0$,
$\mathrm{char}(k)=2$, $K$ has {\it no} places of degree $1$
and $\delta=2$, there are $G_v\cong F^*$, where $|F:k|=4$.
See [6, Lemma 2.5] (with $\mathrm{char}(k)=2$ and $\tau=0$).
See also [6, Remark 2.7].
\par
The case $m=1$ can also occur. Let for
example $k=\FF_2(t)$ and $\CCC=k[X,Y]$ with $Y^2 +XY=X^3 +t$. Then
by the results of Takahashi [12] the part of the quotient graph
corresponding to $X=0$ is a vertex with stabilizer $k(\sqrt{t})$.
\\ \\
\noindent The existence of stabilizers of the type $F^*$, where $F$
is an inseparable extension of $k$ and $|F:k|> 4$ remains an open
question. However the existence of such stabilizers is equivalent to
the existence of a constant field extension of the type described in
the proof of Theorem 2.13 for the following reason. If $L/K$ is a
quadratic extension with constant field $F$, then the stabilizer of
the vertex defined by the lattice class of
$\OOO_{\infty}\oplus\OOO_{\infty}\alpha$ with $\alpha\in F\setminus k$
is $F^*$. It does not seem to be known whether there exist function
fields $K/k$ in characteristic $p$ such that a suitable extension
$L/K$ of degree $p$ has a constant field $F$ of degree $p^m$ over
$k$ with $m>2$.
\\ \\
\noindent {\bf Terminology.} \begin{itemize} \item[(i)] When $G_v$
is determined by Theorem 2.1(a) we call it {\it rational
non-abelian}.
\item[(ii)] When $G_v$ is determined by Theorem 2.3 we call it {\it rational abelian}.
\item[(iii)] When $G_v \cong L^*$, as in Corollary 2.7, we call it of type {\it CM} (for ``complex
multiplication").
\item[(iv)] When $G_v \cong GL_2(k)$, as in Corollary 2.8 and
Theorem 2.1(b), we call it {\it split quaternionic}.
\item[(v)] When $G_v \cong Q^*$, as in Corollary 2.9, we call it
{\it non-split quaternionic}.
\item[(vi)] When $G_v$ is determined by Theorem 2.13 we call it {\it inseparable}.
\end{itemize}
\noindent Note that if $u,v$ are $G$-equivalent they are of the same
type. We conclude this section by separately recording the results
for the case where $k$ is quadratically closed. From all the
previous results of this section (including the proof of Theorem
2.1) the situation here
is considerably simplified.
\\ \\
\noindent {\bf Theorem 2.14.} \it Suppose that $k$ is quadratically
closed. Then either $G_v$ is rational or split quaternionic.
\\ \\
\noindent \rm Theorems 2.1, 2.3, 2.6 and 2.13 apply to every edge stabilizer $G_e$ by virtue of Lemma 1.3. Our results show however that the structure of $G_e$ is restricted in some interesting ways.\\ \\
\noindent {\bf Proposition 2.15.} \it Suppose that $\delta$ is odd. Let the matrix $M$ lie in $G_e$, where $e \in \edge(\TTT)$.
Then the eigenvalues of $M$ lie in $k$.\\ \\

\noindent {\bf Proof.} \rm Let the
characteristic polynomial of $M$ be
$$t^2+\rho t+\tau,$$
where $\rho= -\mathrm{tr}(M)\in k$ and $\tau =\det(M) \in k^*$.
(See Lemma 1.3.) Assume to the contrary that this polynomial is irreducible over $k$. Now there exists $g \in GL_2(K)$ such that $$ gMg^{-1}=\left[\begin{array}{lll} 0 & -1\\[10pt]
\tau & -\rho
\end{array}\right].$$
\noindent Let $C=<M>$. Then $ C'=gCg^{-1} \leq GL_2(k)$ and so $$C' \leq
G_{g(v)} \cap G_{g(v')},$$
\noindent where $v,v'$ are the endpoints of $e$ (in $\TTT$). \noindent  Suppose now that $$\left[\begin{array}{lll} 0 & -1\\[10pt]
\tau & -\rho
\end{array}\right]\in G(n,z),$$
\noindent for some $n,z$. We now apply Lemma 1.2. By parts (i) and (ii) it follows that $\nu(z),n \geq 0$. Suppose now that $n >0$.
Then by Lemma 1.2(iii) it follows that $ z \equiv \lambda\; (\mod \mm)$, for some $\lambda \in k$, since $\delta$ is odd. (Recall that the degree of
$\OOO/\mm$ over $k$ is $\delta$.) This
contradicts the irreducibility of the above polynomial. Hence $n=0$ and so $v(n,z)=v(0,0)$.
It follows that $g(v)=g(v')=v(0,0)$ and hence that $v=v'$. The result follows. \hfill $\Box$ \\ \\
\noindent \rm When $\delta$ is odd therefore the structure of $G_e$ is given by Theorems 2.1 and 2.3. For the arithmetic case Corollaries 2.2, 2.4 and 2.12 ensure that Proposition 2.15
reduces to the following.\\ \\
\noindent {\bf Corollary 2.16.} \it Suppose that $k=\FF_q$ and that $\delta$ is odd.  Then $|G_e|$ is
not divisible by $q^2-1$.\\ \\
\noindent \rm Proposition 2.15 has another interesting consequence. In the above terminology any $G_v$ containing a matrix with eigenvalues not in $k$ must be of types {\it CM},
 {\it quaternionic} or {\it inseparable}. However when $\delta$ is odd $G_v$ cannot be {\it non-split quaternionic} by Lemma 2.10.
 \\ \\
 \noindent {\bf Corollary 2.17.} \it Suppose that $\delta$ is odd and that $G_v$ is CM or inseparable. Then $G_e$ is trivial (i.e. $\cong k^*$),
 for every $e \in \edge(\TTT)$ incident with $v$.\\ \\
 \noindent {\bf Proof.} \rm Theorem 2.13 ensures that in both cases $G_v \cong F^*$, where $F$ is a quadratic extension of $k$. If $G_e$ is not trivial then $G_e=G_v$ by Lemma 1.4, which contradicts Proposition 2.15. \hfill $\Box$ \\ \\
\noindent \rm Nagao's Theorem [9, Corollary, p.87] provides an example (for the case $\delta=1$) of a split quaternionic $G_v$ and an edge $e$ incident with $v$ for which $G_e$ is non-trivial. The restriction on $\delta$
in Proposition 2.15 is necessary.\\ \\
\noindent {\bf Proposition 2.18.} \it Suppose that $k=\FF_q$ and that $\delta$ is even. Let $v \in \vert(\TTT)$ be any vertex for which
$$G_v \cong GL_2(\FF_q).$$
\noindent (e.g. $v=v(0,0)$) Then, for at least one edge $e$ incident with $v$, $|G_e|$ is divisible by $q^2-1$.\\ \\
 \noindent {\bf Proof.} \rm Let $e$ be such an edge. Now $G_v$ acts on the $q^{\delta}+1$ edges incident with $v$ and the order of the
 orbit for this action containing $e$ is $|G_v:G_e|$. (See Lemma 3.1.) If $q^2-1$ does not divide $|G_e|$ then
  $$|G_e|=q^s(q-1)^t,$$
  \noindent where $s\geq 0$ and $t=1,2$ by Corollaries 2.2, 2.4 and 2.6. Then the order of the orbit containing this $e$ is divisible by $q+1$. This cannot apply to all these orbits since then $q^{\delta}+1$ would be divisible by $q+1$, contradicting the fact that $\delta$ is even. The result follows. \hfill $\Box$ \\ \\

\begin{center}{\bf \large 3. The quotient graph: non-isolated vertices}
\end{center}

\space \noindent \rm Throughout this paper the terminology
``degree" has been used in the context of algebraic function fields.
To avoid confusion in the context of graph theory we will use
"valency" instead of "degree" for number of edges attached to a vertex.
\\ \\
\noindent {\bf Notation.} We denote the {\it valency} of a vertex
$v$ in a graph $\mathcal{G}$ by $\val(v)$.
\\ \\
\noindent In this section we are primarily concerned with the values
of $\val(\widetilde{v})$. When $\delta=1$ the situation turns out to be
surprisingly straightforward. In view of Lemma 1.4 is is clear that
all the results of the Theorems 2.1, 2.3, 2.6 and 2.13 apply
to all $G_e$, where $e \in \edge(\TTT)$.
\\ \\
\noindent Now, for each $v\in\vert(\TTT)$, the group $G_v$ acts on,
$\mathrm{star}(v)$, the set of edges of $\TTT$ which are incident
with $v$. Let $\mathrm{Orbs}(v)$ denote the set of equivalence
classes of this action and, for each $e \in \mathrm{star}(v)$, let
$\mathrm{orb}(e)$ be the equivalence class containing $e$.  We recall that
$\TTT$ is a regular graph of valency $\card(k_{\infty})+1$. The fact
that $G$ acts on $\TTT$ {\it without inversion} leads to the following
well-known result.
\\ \\
\noindent {\bf Lemma 3.1.} \it
\begin{itemize}
\item[(i)]For each $v\in \vert(\TTT)$,
$$\mathrm{card}(\mathrm{Orbs}(v)) = \val(\widetilde{v}).$$
\item[(ii)] For each $e \in \mathrm{star}(v)$, there is
a one-one correspondence
$$\mathrm{orb}(e)\; \longleftrightarrow\; G_v/G_e.$$
\end{itemize}\rm
\noindent To repeat from the above the valency of every vertex of
$G\backslash \TTT$ which lies on a ray $R(\sigma)$ (as in Theorem
1.1) is $2$.
\\
\noindent We recall that $\RRRR$ denotes a complete
set of representatives for $\OOO/\mm\;(\cong k_{\infty})$.
\\ \\
\noindent {\bf Lemma 3.2.} \it Let
$$ g=\left[\begin{array}{lll} a & b\\[10pt]
c & d\end{array}\right]\in G(n,z)$$
\noindent and let $u \in \RRRR$. Then $g$ maps
$v(n+1,z+u\pi^n)$ onto
$v(n-1,z)$ if and only if
$$(cz+d)+uc\pi^n \in \mm.$$

\noindent {\bf Proof.} \rm Note that $cz+d,\;c\pi^n \in \OOO$ by
Lemma 1.2. The result follows from Lemma 3.1.
\hfill $\Box$
\\ \\
\noindent Our next result is essential for dealing with rational
stabilizers.
\\ \\
\noindent {\bf Lemma 3.3.} \it Suppose that $G(n,z)$ is rational.
Then there exist $n',z'$ such that
\begin{itemize}
\item[(i)] $v(n',z') \equiv v(n,z)\;(\mod\;G),$
\item[(ii)] $G(n',z')$ consists of all matrices of the form
$M(n',z',\alpha,\beta,c).$
\end{itemize}

\noindent {\bf Proof.} \rm Let $$g=\left[\begin{array}{lll} a & b\\[10pt]
c & d\end{array}\right]\in G(n,z).$$ \noindent From the proofs of
Theorems 2.1 and 2.3 it is clear that every element of $G(n,z)$ is
of the form $M(n,z,-,-,-)$ unless, for all $g \in G(n,z)$, either
(i) $b=0$ or (ii) $c=0$.
\\ \\
\noindent Let $v(n',z')= g_0v(n,z)$, where, for any nonzero $c' \in \CCC$, $$ \mathrm{(i)}\;\; g_0=\left[\begin{array}{lll} 1 & c'\\[10pt]
0 & 1\end{array}\right]\; \mathrm{and\;  (ii)}\;\; g_0=\left[\begin{array}{lll} 1 & 0\\[10pt]
c' & 1\end{array}\right].$$ Then $G(n',z')=g_0G(n,z)g_0^{-1}$ has
the required properties.
\hfill $\Box$
\\ \\
\noindent For now we will focus on rational stabilizers. Our next
result shows that isolated vertices in $G
\backslash \TTT$ never arise from such vertices.
\\ \\
\noindent {\bf Theorem 3.4.} \it Suppose that $G_v$ is rational. Then $$\val(\widetilde{v})\geq 2.$$ .

\noindent {\bf Proof.} \rm Let $v=v(n,z)$. By Lemma 3.3 we may
assume that every element of $G(n,z)$ is of the form
$$M(n,z,\alpha,\beta,c).$$\noindent Then
$$G(n,z) \leq G(n+1,z),$$
\noindent by Lemma 1.2. On the other hand for each
$$g= \left[\begin{array}{lll} a & b\\[10pt]
c & d\end{array}\right] \in G(n,z)$$ \noindent $cz+d \in k^*$. It
follows from Lemma 3.2 that there does not exist $g' \in G(n,z)$
such that $g'v(n-1,z)=v(n+1,z)$. The result follows from Lemma 3.1.
\hfill $\Box$
\\ \\
\noindent We will require the following special case.
\\ \\
\noindent {\bf Lemma 3.5.} \it For any $z$, either $G(0,z) \cong
GL_2(k)$ or
$$G(0,z) = \left\{\left[\begin{array}{lll} \alpha & \beta\\[10pt]
0 &\alpha\end{array}\right]: \alpha \in k^*,\; \beta \in k
\right\}.$$

\noindent {\bf Proof.} \rm If $\nu(z) \geq 0$, then
$G(0,z)=G(0,0)=GL_2(k)$. We may suppose then that $\nu(z)<0$.
\\ \\
\noindent Let $$g=\left[\begin{array}{lll} a & b\\[10pt]
c & d\end{array}\right]\in G(0,z).$$ \noindent Now $ c \in k$ by
Lemma 1.2. Suppose that $ c \neq 0$. By Lemma 1.2 $\nu(a-cz) \geq
0$. It follows that $G(0,z)= G(0,a')$, for some $a' \in \CCC$. Then
$g'v(0,a')= v(0,0)$, where
$$g'=\left[\begin{array}{cc} 1 & -a'\\[10pt]
0 & 1\end{array}\right].$$

\noindent We may assume then that $c =0$. Hence by Lemma 1.2 $\;a,d
\in k^*$ and $\nu(b+z(a-d)) \geq 0$. If $a \neq d$ then by a
previous argument $G(0,z) \cong GL_2(k)$. Otherwise $a=d$ and $ b
\in k$.
\hfill $\Box$
\\ \\
\noindent The latter possibility occurs, for example, when $\nu(z)$
is a {\it gap number} for $\nu$ [11, I.6.7, p.32]. From now on
nearly all our results will apply to the case where $\delta=1$,
i.e. when $\OOO/\mm \cong k$. Here a stronger version of Theorem 3.4 holds.
\\ \\
\noindent {\bf Theorem 3.6.} \it Suppose that $\delta=1$. If $G_v$
is rational non-abelian, then
$$\val(\widetilde{v})=2\;or\;3.$$

\noindent {\bf Proof.} \rm Let $v=v(n,z)$. By Lemma 3.3 we may
suppose that every element of $G(n,z)$ is of the form
$M(n,z,-,-,-)$. If $$g=\left[\begin{array}{lll} a & b\\[10pt]
c & d\end{array}\right]\in G(n,z),$$ \noindent where $n<0$, then
$c=0$ by Lemma 1.2. It follows form this and Lemma 3.5 that $n \geq
1$. As in Theorem 3.4
$$G(n,z) \leq G(n+1,z).$$
\noindent Let $\rho, \epsilon \in k^*$. Choose
$g_0=M(n,z,\alpha_0,\beta_0,-) \in G(n,z)$ such that
$\alpha_0\rho=\beta_0\epsilon$. Then it is easily verified that
$$g_0v(\pi^{n+1}, z+\rho\pi^n)=v(\pi^{n+1},z+\epsilon\pi^n),$$
\noindent i.e.
$$v(\pi^{n+1}, z+\rho\pi^n) \cong v(\pi^{n+1},z+\epsilon\pi^n)\; (\mod\; G(n,z)).$$
\noindent However by Lemma 3.2
$$v(n-1,z) \not\cong v(n+1,z)\;(\mod\; G(n,z)).$$
\noindent The result follows from Lemma 1.8.
\hfill $\Box$
\\ \\
\noindent We repeat that the valency of the projection of a vertex
of $\TTT$ which projects onto those of one of the rays (in Serre's
Theorem) is $2$. Takahashi [12] however has shown that the valency
$3$ can occur for rational non-abelian $G_v$ (when $g=\delta=1$).
\\ \\
\noindent {\bf Theorem 3.7.} \it Suppose that $\delta=1$. If $G_v$
is rational abelian, then
$$\val(\widetilde{v})=2\ \ \ \hbox{\it or}\ \ \
\val(\widetilde{v})=\card(k)+1.$$

\noindent {\bf Proof.} \rm If $G_v$ is trivial (i.e. $G_v \cong
k^*$), then by Lemma 3.1 each $G_v$-orbit on $\mathrm{star}(v)$
contains only one element. We may suppose then by Lemma 3.3 that
$G_v$ consists of elements of the form $M(n,z,\alpha,\beta,c)$,
where for some $c\neq 0$. By an argument used in the proof of
Theorem 3.6,
together with Lemma 3.5, we may further assume that $n \geq 1$.\\ \\
\noindent As before $G(n,z) \leq G(n+1,z)$ and
$$v(n-1,z)\not\equiv v(n+1,z)\;(\mod\;G(n,z)),$$
\noindent by Lemma 3.2. Suppose that there exists $c_0 \in \CCC(n)
\backslash \CCC(n-1)$ for which $M(n,z,\alpha,\beta,c_0) \in
G(n,z)$, where $\alpha \neq \beta$. Then using this element in
conjunction with Lemma 1.4 it follows from Lemma 3.2 that
$$v(n-1,z) \equiv v(n+1,z+\lambda\pi^n),$$
\noindent for all $\lambda \in k$. In this case
$$\val(\widetilde{v})=2.$$
\noindent We are left with the case where $ c \in \CCC(n-1)$, for
all $M(n,z,\alpha,\beta,c)\in G(n,z)$. Here it is easily verified
that $G(n,z) \leq G(n-1,z)$ and that
$$G(n,z) \leq G(n+1,z+\lambda\pi^n),$$
\noindent for all $\lambda \in k$. By Lemmas 3.1 and 1.8 it follows
that $$\val(\widetilde{v})= \card(k)+1.$$
\hfill $\Box$
\\ \\
\noindent Examples of vertices with trivial stabilizers are given in
the previous section. In addition Takahashi [12] (for the case
$g=\delta=1)$ has shown the existence of rational abelian
stabilizers $G_v \;(\cong k^* \times k^+)$ for which
$\deg(\widetilde{v})=2$.
\\ \\

\begin{center}{\bf \large 4. The quotient graph: isolated vertices}
\end{center}

\smallskip \noindent {\bf Definition.}  A vertex $\widetilde{v} \in G
\backslash \TTT$ is {\it isolated} if and only if
$$\val(\widetilde{v})=1.$$
\noindent \rm By Lemma 3.1 for such a vertex we have

$$G_v/G_e \longleftrightarrow \card(k_{\infty})+1,$$
\noindent for all $e \in \mathrm{star}(v)$. Isolated vertices always
exist when $\delta=1$. For example the image of $v(0,0)$ is always
isolated in this case [9, Exercise 6, p.99]. (Its stabilizer is
$GL_2(k)$.) On the other hand $G\backslash \TTT$ need not have any
isolated vertices. If $g=0$ and $K$ has a place of degree $1$ (so
that $K\cong k(t)$) then the results of [4] show that this holds
for all $\delta>1$. We will prove that when $\delta=1$ isolated
vertices arise from non-rational stabilizers.
\\ \\
\noindent {\bf Theorem 4.1.} \it Suppose that $\delta=1$ and that,
for some $ v \in \vert(\TTT)$, $G_v$ contains a matrix $M$ with no
eigenvalues in $k$. Then
\begin{itemize}
\item[(i)]$\val(\widetilde{v})=1$;
\item[(ii)]
$$I(M) \cap G_e=Z,$$
for all $e \in \mathrm{star}(v)$ \end{itemize}

\noindent {\bf Proof.} \rm Let $v=v(n,z)$ as above and let $e$ be
the edge joining $v(n,z)$ and $v(n-1,z)$. Now $G_v$ contains all
matrices of the form
$$\alpha I_2+\beta M,$$ where $\alpha,\beta \in k$, with
$(\alpha,\beta)\neq (0,0)$, by Lemma 1.4. Let
$$M=\left[\begin{array}{lll} a & b\\[10pt]
c & d\end{array}\right].$$

\noindent Let $u \in \RRRR$. Since $\delta =1$ it is then clear that
there exist $\alpha'\in k,\beta'\in k^*$ such that
$$ \alpha'+\beta'(d+cz)+u\beta' \pi^n c \in \mm.$$ \noindent Then
by Lemmas 1.8, 3.1 and 3.2 all edges in $\mathrm{star}(v)$ are
$G(n,z)$-equivalent and so $\widetilde{v}$ is isolated.
\\ \\
\noindent For part (ii) suppose that $I(M) \cap G_e\neq Z$. Then
$I(M) \leq G_{v_0}$, where $v_0=v(n-1,z)$, by Lemma 1.5, which
contradicts the proof of part (i). Hence $I(M) \cap G_e = Z$. Let
$e' \in \mathrm{star}(v)$. From the proof of part (i) there exists
$g \in G_v$ such that $e'=g(e)$. Part (ii) follows.
\hfill $\Box$
\\ \\
\noindent Theorem 4.1 applies to all quaternionic and inseparable
$G_v$, where $k$ is {\it not} quadratically closed. \noindent The
restriction on $\delta$ in Theorem 4.1 is necessary. For the case
where $g=0$, $K$ has no degree $1$ places and $\delta =2$ it is
known [6, Corollary 2.10] that there exists $e \in \edge(\TTT)$
for which $G_e=I(M)\;(\neq Z)$. It is also known from this case
[6, Lemma 2.8] that isolated vertices can occur when $\delta \neq
1$. We do however prove the following partial converse.
\\ \\
\noindent {\bf Theorem 4.2.} \it Suppose that $G_v = I(M)$, where
$M$ has no eigenvalues in $k$. Then
$$\val(\widetilde{v})=1\;\Longleftrightarrow\; \delta=1.$$

\noindent {\bf Proof.} \rm Now $\widetilde{v}$ is isolated when $\delta=1$ by
Theorem 4.1. Assume then that $\widetilde{v}$ is isolated. Let $v=v(n,z)$ and
$$M=\left[\begin{array}{lll} e & f\\[10pt] g & h\end{array}\right].$$
\noindent By a previous argument we may assume that $n \geq 1$. Let
$u \in \RRRR$. As in the proof of Theorem 2.9 there exist $\alpha
\in k,\beta \in k^*$ such that
$$\alpha+\beta(h+gz)+\beta u\pi^ng \in \mm.$$
\noindent We recall from Lemma 1.2 that $\nu(g) \geq -n$. Suppose
that $\nu(g)>-n$. Then $G_v \leq G(n-1,z)$ by Lemmas 1.2 and 1.4
which implies that $\widetilde{v}$ cannot be isolated. Hence $\pi^ng \in
\OOO\backslash \mm$. It follows that there exist $\gamma \in k$ and
constants $r,s \in \OOO$ such that
$$u \equiv  \gamma r+s\;(\mod \mm).$$
\noindent Hence $\delta$, the degree of $\OOO/\mm$ over $k$, is $1$.
\hfill $\Box$
\\ \\
\noindent For the case where $g=0$, $K$ has no degree $1$ places and
$\delta =2$ it is known there exists $v \in \vert(\TTT)$ for which
(i) $\widetilde{v}$ is isolated and (ii) $G_v \neq I(M)$.
See [6, Lemmas 2.5, 2.8].
\\ \\
\noindent {\bf Corollary 4.3.} \it Suppose that $\delta=1$ and let
$\{ \widetilde{v_\lambda}: \lambda \in \Lambda\}$ be the set of all
$\widetilde{v_\lambda} \in \vert(G \backslash\TTT)$ such that
$$G_{v_\lambda}=I(M_\lambda),$$
\noindent for some $M_\lambda$ with no eigenvalues in $k$.
\noindent
Let $L_\lambda$ be the quadratic extension of $k$ generated by the
eigenvalues of $M_\lambda$. Then there exists an epimorphism
$$\theta: G \twoheadrightarrow \displaystyle{\Star_{\lambda\in\Lambda}}(L_\lambda^*/k^*).$$

\noindent {\bf Proof.} \rm  By Theorem 4.1
$$\val(\widetilde{v_\lambda})=1.$$
Let $e_\lambda \in\edge(\TTT)$ be a lift of its adjacent edge to one
adjacent to $v_\lambda$ ($\in \vert(\TTT)$).  Then $G_{e_\lambda}=Z$ by
Theorem 4.1.
\\ \\
\noindent From the fundamental theorem of the theory of groups
acting on trees [9, Theorem 13, p.55] $G$ is isomorphic to the {\it
fundamental group of the graph of groups} given by $G \backslash
\TTT$. (See [9, p.42].) Hence, for each $\lambda \in \Lambda$,
$$G /Z \cong (G_{v_\lambda} /G_{e_\lambda}) \star H,$$
where $H$ is non-trivial. The result follows.
\hfill $\Box$
\\ \\
\noindent Even when $\delta=1$ it can happen that
$\Lambda=\emptyset$. For the simplest case consider $g=0,\;\delta=1$
(so that $\CCC \cong k[t]$). By {\it Nagao's theorem}, [9,
Corollary, p.87] if $G_v$ is
non-rational then, $G_v$ is split quaternionic.\\ \\
\noindent However Corollary 4.3 does occur. Suppose, once again,
that $g=\delta=1$. Then $\CCC$ is the coordinate ring of an {\it
elliptic} curve. In which case there exist $X,Y \in \CCC$, where
$\nu(X)=-2$ and $\nu(Y)=-3$ which satisfy a {\it Weierstrass
equation}, $F(X,Y)=0$, for which
$$\CCC=k[X,Y].$$\noindent Takahashi [12, Theorem 5]
has shown that, if for some $\lambda\in k$ there does {\it not}
exist any $\mu \in k$ for which $F(\lambda,\mu)=0$, then
$$G(2,\pi^{-1}+ \pi\lambda) \cong k(\omega)^*,$$ \noindent
where $k(\omega)$ is the quadratic extension of $k$ given by any
$\omega$ for which $F(\lambda,\omega)=0$. An explicit example can be
found in [12, p.87].\\ \\

\begin{center}
{\bf \large 5. Special constant fields}
\end{center}

\noindent We list separately the results for the arithmetic case.
\\ \\
\noindent {\bf Theorem 5.1.} \it Suppose that $k=\FF_q$. Let $ v \in
\vert(\TTT)$. Then $\val(\widetilde{v})=1$ if and only if
\begin{itemize}
\item[(i)] $\delta=1$,
\item[(ii)] $\FF_{q^2}^* \hookrightarrow G_v.$
\end{itemize}
\noindent {\bf Proof.} \rm The conditions are sufficient by Theorem
4.1 and Corollary 2.12. Suppose now that $\widetilde{v}$ is isolated. Then
$q^{\delta}+1$ divides $|G_v|$ by Lemma 3.1. (We recall that
$\OOO/\mm \cong k_{\infty}$ and that $|k_{\infty}:k|=\delta$.) From the
list of all possible $|G_v|$ listed in Corollaries 2.2, 2.4 and 2.12
it follows that $\delta=1$. If the second condition is not satisfied
then by the same results $G_v$ is of rational type. Then
$\widetilde{v}$ cannot be isolated by Theorem 3.4.
\hfill $\Box$
\\ \\
\noindent The significant difference here from the non-arithmetic case is
that when $k$ is finite isolated vertices arise {\it only} when
$\delta =1$. We repeat that, from [6, Lemma 2.8], it is possible,
when $k$ is infinite, for $G \backslash \TTT$ to have an isolated
vertex when $\delta > 1$. \noindent From the results of the two
previous sections we can state.
\\ \\
\noindent {\bf Corollary 5.2.}
\it Suppose that $k=\FF_q$ and that $\delta=1$. Let $v \in
\vert(\TTT)$. Then $$\val(\widetilde{v}) \;=\;1, 2, 3 \;or\; q+1.$$

\noindent \rm The possible patterns of orbit sizes are,
respectively, $(q+1)$, $(1,q)$, $(1,1,q-1)$ and $(1,\cdots,1)$.
Takahashi [12] has shown (for the case $g=\delta=1$) that all
possibilities can occur. We now record the arithmetic version of
Corollary 4.3. Recall that here $X$ is {\it finite}.
\\ \\
\noindent {\bf Corollary 5.3.} \it Suppose that $k=\FF_q$ and that
$\delta=1$. Let $\{ \widetilde{v_i}: 1 \leq i \leq n\}$ be the set
of all $\widetilde{v_i} \in \vert(G \backslash\TTT)$ such that
$$G_{v_i}\cong \FF_{q^2}^*.$$
\noindent Then there exists an epimorphism $$\theta: G
\twoheadrightarrow \ZZ/(q+1)\ZZ \star\cdots \star\ZZ/(q+1 )\ZZ,$$
to $n$ factors.
\\ \\
\noindent \rm  Takahashi [12, p.87] provides an example of Corollary
5.3, for the case $k=\FF_3$, where $n=1$. Here then
$$\theta:G/Z \twoheadrightarrow \ZZ/4\ZZ.$$
\noindent In [7] we will, among other things, derive more
information about $n$. This requires, however, quite a lengthy
detour and the use of arguments from algebraic number theory.\\ \\
\noindent In conclusion we deal with the one case yet to be
considered.
\\ \\
\noindent {\bf Theorem 5.4.} \it Suppose that $k$ is quadratically
closed and that $\delta=1$. If $G_v$ is split quaternionic, (i.e.
$\cong GL_2(k)$), then
$$\val(\widetilde{v})=1.$$

\noindent {\bf Proof.} \rm From the proof of Theorem 2.1, part (b),
by conjugating $G'_v$ by $$\left[\begin{array}{cc} 0 & b_0\\[10pt]
c_0 & 0\end{array}\right],$$ (where $b_0c_0=1$) it follows that
$G_v$ is conjugate in $GL_2(K)$ to $GL_2(k)$. Then
$$GL_2(k)=G(n,z),$$
for some $n,z$.
\\ \\
\noindent  As before we may assume that $n \geq 0$ since $GL_2(k)$
is not metabelian. It follows from Lemma 1.2 that $\nu(z) \geq 0$.
Suppose that $n \geq 1$. Then by Lemma 1.2 $$ b+z(a-d)-z^2c \in
\mm,$$ \noindent for all $$\left[\begin{array}{lll} a & b\\[10pt]
c & d\end{array}\right]\in GL_2(k).$$ \noindent Hence $n=0$ and so
$$G(n,z)=G(0,0).$$
\noindent The result follows from [9, Exercise 6, p.99].
\hfill $\Box$
\\ \\

\begin{center}
{\bf \large 6. Concluding remarks}
\end{center}

\noindent Let $\TTT_m$ be a maximal subtree of $G\backslash \TTT$
and let $$\omega=\mathrm{card}(\{e \in \edge(G \backslash \TTT): e
\notin \edge(\TTT_m)\})=\mathrm{card}(\{e \in \edge(X): e \notin
\edge(\TTT_m)\}).$$ \noindent Now we put
$$G_V=\langle G_v: v \in \vert(\TTT)\rangle.$$

\noindent From standard Bass-Serre theory [9, p.43] it is known
that
$$G/G_V \cong F_{\omega},$$
the free group of rank $\omega$. This raises the following question.
\begin{center}
\noindent {\it When is $G \backslash \TTT$ a tree ?}
\end{center}

\noindent The theory of Drinfeld modular curves provides a complete
answer [5] when $k$ is finite. When $k$ is infinite very little is
known. If $g=0$ and $K$ has a place of degree one (so that $K\cong k(t)$)
a complete answer is known [4]. It is also known [6] that
$G\backslash \TTT$ is a tree when $g=0$, $K$ has no places of degree
one (i.e. $K \not\cong k(t)$) and $\delta=2$. Takahashi [12] has
shown that the quotient graph is a tree for the {\it elliptic} case
$g=\delta=1$.

\subsection*{\hspace*{10.5em} References}
\begin{itemize}
\item[{[1]}] V.~G.~Drinfeld, Elliptic Modules,  Math. USSR-Sbornik
 23  (1976) 561-592.
\item[{[2]}] E.-U.~Gekeler,  Drinfeld Modular Curves,
Lecture Notes in Mathematics vol. 1231 Springer, Berlin, Heidelberg,
New York, 1986.
\item[{[3]}] A.~W.~Mason, Serre's generalization of Nagao's
theorem: an elementary approach,  Trans. Amer. Math. Soc.
 353  (2003) 749-767.
\item[{[4]}] A.~W.~Mason, The generalization of Nagao's
theorem to other subrings of the rational function field,  Comm.
Algebra 31 (2003) 5199-5242.
\item[{[5]}] A.~W.~Mason and A.~Schweizer, The minimum index
of a non-congruence subgroup of $SL_2$ over an arithmetic domain,
 Israel J. Math.  133  (2003) 29-44.
\item[{[6]}] A.~W.~Mason and A.~Schweizer, Nonrational Genus
Zero Function Fields and the Bruhat-Tits Tree,  Comm. Algebra
37  (2009) 4241-4258.
\item[{[7]}] A.~W.~Mason and A.~Schweizer, Elliptic points
of the Drinfeld modular groups,  in preparation.
\item[{[8]}] R.~Pierce,  Associative Algebras, Springer
GTM 88, Berlin, New York, 1982.
\item[{[9]}] J.-P.~Serre,  Trees,  Springer, Berlin,
Heidelberg, New York, 1980.
\item[{[10]}] H.~Stichtenoth, \"Uber das Geschlecht eines
inseparablen Funktionenk\"orpers,  Manuscripta Math.  14
 (1974) 173-182.
\item[{[11]}] H.~Stichtenoth,  Algebraic Function Fields and
Codes,  Springer, Berlin, Heidelberg, New York, 1993.
\item[{[12]}]S.~Takahashi, The fundamental domain of the
tree of $GL(2)$ over the function field of an elliptic curve,
 Duke Math. J.  72  (1993) 85-97.
\item[{[13]}] G.~Villa Salvador,  Topics in the Theory of Algebraic
Function Fields,  Birkh\"auser, Boston, 2006.
\end{itemize}

\end{document}